\documentclass{article}
\usepackage{amssymb}

\begin {document}

LANL e-print math.AP/0203230

\bigskip

\centerline {\bf Olga S. Rozanova\footnote{Partially supported by
RFBR grant 00-02-16337}}

\centerline {\it Moscow State University}

\smallskip

\begin {center}
{\large \bf
On classes of globally smooth
solutions to the Euler equations in several dimensions
}
\end {center}
\newtheorem{theo1}{Theorem 1.}
\newtheorem{theo2}{Theorem 2.}
\newtheorem{cor}{Corollary 2.}

\bigskip
The aim of the present work is the investigation of solutions to
the system of the Euler equations, including the systems with the
right-hand sides describing various interior forces:
$$ \partial_t \rho + {\rm div}\, (\rho{\bf V})=0, \eqno (1)
$$ $$ \partial_t (\rho{\bf V})+(\rho{\bf V},{\bf \nabla})\,{\bf
V}+{\bf\nabla} p = \rho {\bf f}({\bf x}, t, {\bf V}, \rho, S), \eqno
(2)$$ $$ \qquad \partial_t S +({\bf V},{\bf\nabla} S)=0 \eqno (3)$$
with the state equation
$p = e^S \rho^\gamma,\,\gamma=const>1.$ We suppose that ${\bf
f}$ is the smooth function of all its arguments.
Here $\rho(t,{\bf x}) ,\,{\bf
V}(t,{\bf x}),\, S(t,{\bf x})$ are the components of solution,
corresponding to the density, the velocity and the entropy,
given in
${\mathbb R}\times {\mathbb
R}^n,\, n\ge 1.$
Equations (1--3) describe the balance of mass, momentum and entropy,
correspondingly.

Set the Cauchy problem for (1 -- 3):

$$ \rho(0,{\bf x})=\rho_0({\bf x})\ge 0,
\, {\bf V}(0,{\bf x})={\bf V}_0({\bf x}),\, S(0, {\bf x})=S_0({\bf
x}).\eqno (4)$$

We deal with the classical solutions to system (1--3)
with the density so quickly decreasing as  $|{\bf x}|\to \infty,$
to guarantee the convergence of the integral
$\int \limits_{{\mathbb R}^n}\rho |{\bf x}|^2 d{\bf
x}$ (so called solutions with the finite momentum of inertia).

As well known, the solution to the Cauchy problem for system
(1-3) may lose the initial smoothness in a finite time,
sometimes there is a possibility to estimate the time
of singularity formation from above
(see,
f.e., \cite{RozDU} and references therein).
Moreover, in the case more frequently investigated
${\bf f}=0$ the compactly supported initial data are sufficient
for the further singularity formation.
(f.e., \cite {Xin}).

At the same time it is interesting that there are
some nontrivial classes of globally smooth solutions.
Note that if ${\bf f}=0,$ the components of
such solutions do not belong to the Sobolev class.

The paper is organized as follows. Firstly  {\it supposing} the
existence of globally in time smooth solution having
concrete properties (denote it $U_0$), we show that if we choose
some initial data ($CD_1$) close in the  Sobolev norm to initial data
of such solution ($CD_0$), then the corresponding solution to the
Cauchy problem ($U_1$) occurs globally smooth as well.  Then we show
that the solutions $U_0$ with such kind of properties {\it exists},
moreover, at $n=2$ we construct some of them.

The present work may be considered as a continuation of papers
 \cite{Serre},
\cite{Grassin}, and their generalization in some sense.

\section
{Symmetrization and the result on a local-in-time smoothness}

The result on a local in time existence of smooth
solution
${\bf U}(t,{\bf x})$ for symmetric hyperbolic systems,
i.e.
systems of the form
$$ A^0(t, {\bf x}, {\bf
U})\frac{\partial {\bf U}}{\partial t}+\sum\limits^n_{k=1} A^k(t,
{\bf x}, {\bf U})\frac{\partial {\bf U}}{\partial {\bf x}^k}=g(t,{\bf
x}, {\bf U}),\eqno(1.1) $$
where the matrices $A^j(t,{\bf x}, {\bf U})$
are symmetric, and additionally the matrix $A^0(t,{\bf x}, {\bf U})$
is positive definite,  is well known
(\cite{Kato},\cite{Volpert}, the particular cases in
\cite{GardingLeray}, \cite{Lax}, \cite{Majda})).

Namely, if the matrices $A^j(t,{\bf x},
{\bf p}), g(t,{\bf x}, {\bf p})$ depending smoothly on their
arguments,  have continuous and bounded derivatives
with respect to the variables $({\bf x},
{\bf  p})$
up to order $m+1$ under bounded $\bf p$, the
initial data ${\bf U}_0={\bf U}(0,{\bf x})$ and the function
$g(t,{\bf x}, {\bf 0})$ belong to the Sobolev class
$H^m({\mathbb R}^n),\,m>1+n/2,$ for any fixed $t\ge 0$, then locally
in time the corresponding Cauchy problem has a unique solution from
the class $\cap _{j=0}^1 C^j([0,T);H^{m-j}({\mathbb R}^n)), \,T>0.$
Moreover, $$\lim\limits_{t\to T-0}\sup(\|{\bf U}\|_{L^\infty}+
\|\nabla_{\bf x} {\bf U}\|_{L^\infty})=+\infty,\eqno(1.2)$$ $T=T({\bf
U}_0) $ and $ \displaystyle\lim\limits_{\|{\bf U}_0 \|_{H^m}\to
 0}T({\bf U}_0)=+\infty.$

(Note that  less
rigid requirements may be imposed on the initial data
for the existence of smooth solution locally in time,
f.e. local (in space), but uniform
belonging to $H^m({\mathbb R}^n).$
\cite{Volpert},\cite{Chemin}).

In many problems having the origin in physics
the coefficients of system (1.1) depend only on solution,
this simplifies significantly the formulation of the result.
In the case for the local in time existence of the Cauchy problem
in the Sobolev class described above, it is sufficient, beside of
the coefficients smoothness, to require the implementation of the
condition $g({\bf 0})={\bf 0}.$

As we deals with the solutions to system (1--3) such that
$\inf \rho=\inf p= 0,$ we have to use the symmetrizartion
firstly proposed in
\cite{MakinoUkaiKawashima}.

Involving the variable $\displaystyle
\Pi=\kappa\left(p/2\right)^{\frac{\gamma-1}{2\gamma}},\,
\kappa=\frac{2\sqrt{\gamma}}{\gamma-1},$
we get the symmetric form as follows:

$$ \exp(\frac{S}{\gamma
})(\partial_t +{\bf V},{\bf \nabla})\Pi + \frac {\gamma
-1}{2}\exp(\frac{S}{\gamma})\Pi div\, {\bf V}=0,\eqno(1.3)
$$
$$ (\partial_t +({\bf V},{\bf \nabla})){\bf V}+\frac
{\gamma -1}{2}\exp (\frac{S}{\gamma})\Pi{\bf\nabla} \Pi = {\bf
f}_1( t, {\bf x},\Pi, {\bf V}, S)=
{\bf f}( t, {\bf
x},e^{-\frac{S}{\gamma}}\left(\frac{\Pi}{\kappa}\right)^{\frac{2}{\gamma-1}},
{\bf V}, S),\eqno(1.4)
$$
$$
(\partial_t +({\bf V},{\bf \nabla})) S=0.\eqno(1.5) $$

Denote ${\bf U}= (\Pi,\,{\bf V},\,S).$
As an immediate corollary of the general result on the
symmetric hyperbolic systems we obtain the following

\begin{theo1}
Let the initial data ${\bf U}_0=
(\Pi_0,\, {\bf V}_0,\, S_0),$
belong to the class $H^m({\mathbb R}^n),\,m>1+n/2.$
Suppose the function
${\bf
f}_1( t, {\bf x},{\bf U}))$
have continuous and bounded derivatives
with respect to space variables
and the solution components up to order $m+1$ for bounded
$\bf U,$
and
$
{\bf
f}_1( t, {\bf x},{\bf 0})
\in H^m({\mathbb
R}^n)$  for any fixed $t\ge 0$.

Then the Cauchy problem for system (1 --
3) has locally in time a unique solution such that
$$(\Pi,\,
{\bf V},\, S,) \in \cap _{j=0}^1 C^j([0,T);H^{m-j}({\mathbb R}^n),$$
moreover $$ \lim_{\|{\bf U}_0 \|_{H^m}\to 0}T=+\infty.$$
\end{theo1}

\section{The interior solution}

We call {\it an interior solution} globally in time smooth solution
$(\bar
\rho(t,{\bf x}),\, \bar{\bf V}(t,{\bf x}), \,\bar S(t,{\bf x}))$
to system (1--3) having the following property: another
solution $(\rho,\, {\bf V},\,S)$ to the Cauchy problem (4) having
sufficiently small
norm $$\|(\rho_0^{(\gamma-1)/2}({\bf x})- \bar
\rho^{(\gamma-1)/2}(0,{\bf x}),\,{\bf V}_0({\bf x})-\bar {\bf V}(0,
{\bf x}),\, S_0({\bf x})-\bar S(0,{\bf x}) \|_{H^m({\mathbb R}^n)}$$
is smooth globally in time  as well, moreover
$$(\rho^{(\gamma-1)/2}- \bar \rho^{(\gamma-1)/2}, \, {\bf V}-\bar {\bf
V}, \,S-\bar S )\in \cap _{j=0}^1 C^j([0,\infty);H^{m-j}({\mathbb
R}^n)),\, m>1+n/2.$$

Note that the trivial solution is not interior
if
$\bf f
=0$.

The set of interior solutions is not empty.
In the paper \cite{Grassin} (the generalization of
 \cite{Serre}) for ${\bf f}={\bf 0}$ it is shown that
the solution to system (1--3) $(0,\, \bar {\bf
V}(t,{\bf x}),\,const)$ is  interior, if  $\bar {\bf V}(t,{\bf
x})$ is the solution to equation $\partial_t {\bf V}+({\bf V},{\bf
\nabla})\,{\bf V} = 0$ such that $Spec D\bar{\bf
V}(0,{\bf x})$ is separated from  the semi-axis ${\mathbb R}_-$ and
$D\bar{\bf V}(0,{\bf x})\in L^\infty({\mathbb R}^n),\, D^2\bar{\bf
V}(0,{\bf x})\in H^{m-1}({\mathbb R}^n).$ (We denote $D^k$ the vector
of all spatial derivatives of order $k$.)

In particular, as it is noted in \cite{Serre} in the role of
such
solution may be played by the solution with linear profile of velocity
${\bf V}=A(t){\bf
r}, $ where ${\bf r}$ is the radius-vector of point, $A(t)
$ is the matrix with the coefficients depending on time
such that ${\rm Spec } A(0)\notin
{\mathbb R}_-.$ The solution may be found explicitly.

Note that the given interior solution
(and the close solution as well)
does not belong to the Sobolev space,
since the velocity components grow at infinity.
The result is clear from the physical point of view:
the velocity field having the positive divergency
"spreads out" the initially concentrated
sufficiently small mass, that prevent the singularity
formation.

But  the given solution with the zero density is not physical.
It follows, in particular,
in the case of the velocity Jacobian having initially negative spectrum
at some value of $x$ , the appearance during a finite time of
points with infinite negative divergency.
If the solution with non-zero velocity, initially close to the
described above, would behave like this, it would follow
(according to (1)) the infinite value of density in some
points.

The question arises: can one construct the interior solution
with the density not close to zero in the Sobolev norm?  Or to be
interior the solution must describe the scattering of sufficiently
small mass?

If from some reasons we succeed to find the velocity field,
then from the linear with respect to
$\rho$ and $S$ equations (1), (3) we can find the solution
to the Cauchy problem (4), moreover, if the velocity field
occurs smooth, then the rest of the components
are smooth.

It needs to note, that the values
$\rho$ and $p$ can become unbounded
in some point of trajectory in a finite time $T$ if and only if
$\int_0^T div {\bf V}
dt=-\infty$ on the trajectory (see (1)).

Before formulating the theorem we do some transformation
under supposition that
$(\bar\Pi,\,\bar{\bf V}=A(t){\bf r},\,\bar S)$ is a globally in time
smooth solution to system
(1.3--1.5).

According to (1.3--1.5) the vector-function ${\bf
u}={\bf U}-\bar{\bf U}$ (here ${\bf u}=(\pi,\,{\bf
v},\,s),\,$${\bf U}= (\Pi,\,{\bf V},\,S),\,$$\bar{\bf U}=
(\bar\Pi,\,\bar{\bf V},\, \bar S)$) satisfies the system
of equations
$$
(\partial_t +\bar{\bf V},{\bf \nabla})\pi +({\bf v},{\bf \nabla})\pi
+
\frac {\gamma-1}{2}(\pi+\bar\Pi) {\rm div} {\bf V}
=
$$
$$
-({\bf v},{\bf \nabla})\bar\Pi
-\frac {\gamma-1}{2}\pi {\rm div}\, \bar{\bf V},\eqno(2.1)
$$
$$ (\partial_t +(\bar{\bf V},{\bf \nabla})){\bf v}+
({\bf v},{\bf \nabla}){\bf v}+
\frac
{\gamma -1}{2}\exp (\frac{S}{\gamma})(\pi+\bar\Pi){\bf\nabla} \pi =
$$
$$
-({\bf v},{\bf \nabla})\bar{\bf V}
-\frac
{\gamma -1}{2}\exp (\frac{S}{\gamma})(\pi+\bar\Pi){\bf\nabla} \bar
\Pi +
\frac
{\gamma -1}{2}\exp (\frac{\bar S}{\gamma})\bar\Pi{\bf\nabla} \bar
\Pi + \eqno(2.2)
$$
$$+{\bf
f}_1( t, {\bf x}, (\bar{\bf V}+{\bf v}), (\bar S+s),(\bar\Pi+ \pi))-
{\bf
f}_1( t, {\bf x}, \bar{\bf V}, \bar S, \bar\Pi),$$
$$
(\partial_t +\bar{\bf V},{\bf \nabla}) s+
({\bf v},{\bf \nabla}) s
=-({\bf v},{\bf \nabla})\bar S.\eqno(2.3) $$

Further, following to \cite{Serre}, we carry out the nondegenerate
change of variables such that the infinite semi-axis of time
turns out to semi-interval $[0,\sigma_\infty)$,
to obtain the symmetric hyperbolic system with the coefficients
allowing to apply the theorem on a local existence
of smooth solution on the interval
$[0,\sigma_*),$\, $\sigma_*\le\sigma_\infty.$

So, let $(t_1,{\bf y}):= (t, A(t){\bf x})$ be the new variables
($A(t)= \|a_{ij}(t)\|$ is a quadratic $(n\times n)$
matrix).  Then ${\bf\nabla}_{\bf x}=A^*{\bf\nabla}_{\bf
y},\,{\rm div}_{\bf x}{\bf V}= {\rm div}_{\bf y} A{\bf V},
\,\partial_{t_1}=\partial_t+A{\bf r}{\bf\nabla}_{\bf y}$.
Choose the nonnegative function $
\lambda(t)$ such that the integral
$\int_0^\infty \lambda(\tau) d\tau=\sigma_\infty<\infty$ converges
and set
$\sigma(t_1)=\int_0^{t_1} \lambda(\tau) d\tau.$ In that way,
the semi-infinite axis of time goes to the semi-interval
$[0,\sigma_\infty).$

Further, involve the variables
$$
{\bf W}=A(t){\bf v}\lambda^{-1}(t),\,
P=\lambda^{-q}(t)\pi, \, $$
the constant $q$ will be defined below.

After the transformation we get the system
$$
(\partial_\sigma +({\bf W},{\bf \nabla_y}))P +
\frac {\gamma-1}{2}(P+\bar P) {\rm div}_y\, {\bf W}
=
$$
$$
-({\bf W},{\bf \nabla}_y)\bar P
-
PQ_1(\sigma),\eqno(2.4)
$$
$$ (\partial_\sigma +({\bf W},{\bf \nabla}_y)){\bf W}+
\frac
{\gamma -1}{2}\Psi(S,\sigma)
(P+\bar P){\bf\nabla}_y P = $$
$$ -
\frac {\gamma -1}{2}
\Psi(S,\sigma)
(P+\bar P){\bf\nabla}_y \bar P
+
\frac {\gamma -1}{2}
\Psi(\bar S,\sigma)
\bar P{\bf\nabla}_y \bar P
-{\bf
W}Q_2(\sigma)+G,\eqno(2.5) $$
$$ (\partial_\sigma +{\bf W},{\bf \nabla}_y) s =-({\bf
W},{\bf \nabla}_y)\bar S.\eqno(2.6) $$

Here we denote
$$
Q_1(t(\sigma))=\lambda(t)^{-1}(\frac{\gamma-1}{2} {\rm tr} A(t)+
q (\ln \lambda(t))'),$$
$$Q_2(t(\sigma))=\lambda(t)^{-1}((\ln \lambda(t))' E
+A(t)((A^{-1})'(t)+E)),$$
$$\Psi(S,\sigma)=
\exp
(\frac{S}{\gamma})R(\sigma)B(\sigma),$$
$$B(t(\sigma))=
A(t)A^*(t)({\rm det }A(t))^{-2/n}$$
$$R(t(\sigma))= \lambda^{2q-2}(t) ({\rm det}\,A(t))^{2/n},$$
$$G=\lambda^{-2}(t)A(t)({\bf f}_1( t,
A(t){\bf x}, \lambda^{-q}(t)(\bar P+ P),\lambda(t)A^{-1}(t)
(\bar{\bf W}+{\bf W}), (\bar S+s))-$$
$$ {\bf f}_1(
t, A(t){\bf x}, \lambda^{-q}(t)\bar P,
\lambda(t)A^{-1}(t)\bar{\bf W}, \bar S)),$$
taking into account that
$t=t(\sigma).$

Note that  $B$ is a bounded invertible matrix.

After multiplying (2.5) by $\Psi^{-1}(S,\sigma)$
the system become symmetric hyperbolic, however, generally
speaking, $\Psi^{-1}(S,\sigma)\to 0$ as $t\to \infty$
or will be unbounded owing to the behaviour of
$R(\sigma(t))$.

In another words, we can increase the time of existence of
smooth solution to system (1--3) as long as we wish due to choosing
the initial data small in the Sobolev $H^m$ norm, but we cannot
guarantee the existence of the solution during the infinite time.
Let
$$D^\alpha=\left(\frac{\partial}{\partial y_1}\right)^{\alpha_1}...
\left(\frac{\partial}{\partial y_n}\right)^{\alpha_n},\,
\sum\limits_i \alpha_{i=1}^n=\alpha.$$

Denote $C_b^r({\mathbb R}^n)$ the space of continuous
functions having continuous bounded in ${\mathbb R}^n$
derivatives up to order $r=0,1,...$.
Further, denote ${\bf f}_{\bf V}$ the matrix
$\|\frac{\partial f_i}{\partial V_j}\|$.

Now formulate the theorem taking into account all
denotation involved.

\begin{theo2}
Let the function
${\bf
f}_1( t, {\bf x},\Pi, {\bf V}, S)=
{\bf f}( t, {\bf
x},e^{-\frac{S}{\gamma}}\left(\frac{\Pi}{\kappa}\right)^{\frac{2}{\gamma-1}},
{\bf V}, S)$
have the derivatives with respect to all components of
solution up to order
$m+1$, continuous and bounded  under bounded
$(\Pi, {\bf V}, S),$ moreover $f_1(t,{\bf x}, 0, {\bf 0}, 0)\in
H^m({\mathbb R}^n)$ at any fixed $t\ge 0.$
Suppose
system (1--3) has the globally in time smooth solution $\bar
U=(\bar\rho, \bar{\bf V}, \bar S)$ with the linear profile of
velocity $\bar{\bf V}=A(t){\bf r}$ such that

a) $\bar\rho^{\frac{\gamma-1}{2}}(0,{\bf x})\in
\cap
_{j=0}^{m+1} C^j_b({\mathbb R}^n);$
 $D\,\bar S(0,{\bf x})\in
\cap
_{j=0}^m C^j_b({\mathbb R}^n);$

b) $\xi(t)=\det A(t)>0$ for $t\ge t_0 \ge 0;$

\noindent
and there exist a smooth real-valued function
$\lambda(t)$, a constant $q$ and
a skew symmetric matrix with the real-valued coefficients
$U_\phi(t)$
with the following properties:

c)$\int_{t_0}^{+\infty} \lambda(\tau) d\tau<\infty,$

d) $\int_{t_0}^{+\infty} \lambda^q(\tau) \xi^{1/n}(\tau)\,
d\tau<\infty,$

e) functions

$$Q_1(t)R(t),\eqno(2.7)$$

$$Q_2(t)-\lambda^{-1}(A(t){\bf f}_{\bf V}
A^{-1}(t)- U_\phi(t)),\eqno(2.8)$$

$$(\ln R(t))'\eqno(2.9)$$

are bounded on ${\mathbb R}_+\times {\mathbb R}^n$.

Then the solution $\bar U$ is interior.
\end{theo2}

Proof.
At first we consider the case of $R(\sigma)$
such that $0<r_1=const\le R(\sigma)\le r_2=const<\infty$.
If the other conditions of
Theorem 2.1 hold, then we can prove  Theorem 1.1
immediately.  The case is realized, for example, at
$\gamma=\frac{2}{n}+1,\,{\bf f}=0,$ in the situation, considered in
\cite{Serre} ($tr A(t)\sim \frac{n}{t},\, \xi(t)\sim
t^{-n}\,(t\to\infty),\,\lambda(t)=(1+t)^{-2},\,q=\frac{n(\gamma-1)}{4}).$

\vskip1cm
In the general case the proof is more complicated.

Suppose $\sigma\in [0,\sigma_*)$ and $(X,\,{\bf Y},\,Z): R^n\to
R\times R^n \times R $ is a vector-function from $L_2.$
Following \cite{Serre}, involve the norm $$ [X,{\bf
Y},Z]^2(\sigma):=\int\limits_{R^n}(X^2+Z^2+ {\bf
Y}^*\Psi^{-1}(S,\sigma){\bf Y}) d{\bf y}.$$

For $p\in N$ define
$$E_p(\sigma)=\frac{1}{2}\sum\limits_{|\alpha|=p}
[D^\alpha P,\,D^\alpha {\bf W},\,
D^\alpha s]^2(\sigma),$$
$$F_m(\sigma)=\sum\limits_{p=0}^m E_p(\sigma).$$

If we should succeed to show that
$F_m(\sigma) $ is bounded on $[0,\sigma_\infty]$ function,
the theorem would be proved, as the impossibility to prolong
the solution is connected with the going to the infinity of
$L_\infty -$ norms of the solution itself or its gradient (see
(1.2)).

Showing that $H^m -$ norm of the solution admits
a majorization by means of some function, bounded on
$[0,\sigma_\infty]$, we shall show that $\sigma_*=\sigma_\infty$
and thus the theorem will be proved.

After involving of the new variable
$\Pi,$ which is necessary  for the symmetriization of the system
(see above), globally smooth solution
$\bar U$ corresponds to the globally smooth solution
$(\bar\Pi,\,\bar {\bf V},\,\bar S)$
with the initial data
$(\bar\Pi_0,\,\bar {\bf V}_0,\,\bar S_0).$ After the
changing of variables  it  corresponds to the solution
$(\bar P,\,\bar {\bf W},\,\bar S)$ with initial data
$(\bar P_0,\,\bar {\bf W}_0,\,\bar S_0).$ The last solution is
a smooth solution to system
(2.4--2.6) on $[0,
\sigma_\infty]$, and therefore its  $H^m -$ norm is bounded
on the segment.

\vskip1cm
Compute
$$\frac{dE_q}{d\sigma}=
\sum\limits_{|\alpha|=q}[
\int\limits_{R_n}(D^{\alpha}
P D^\alpha\partial_\sigma Pd{\bf y} +$$
$$
\int\limits_{R_n}D^\alpha s D^\alpha\partial_\sigma s d{\bf y}+$$
$$
\int\limits_{R_n}D^\alpha {\bf W} \Psi^{-1}(S,\sigma) D^\alpha\partial_\sigma {\bf W} d{\bf y}+
$$
$$
\int\limits_{R_n}(D^\alpha {\bf W})^*
(\frac{\partial\Psi^{-1}(S,\sigma)}{\partial\sigma}- \frac{1}{\gamma}
\Psi^{-1}(S,\sigma)\partial_\sigma S) D^\alpha {\bf W} d{\bf y}]= $$ $$
I_1+I_2+ I_3 +I_4$$
(further it is supposed the summation over repeated indices,
the integrals $I_k,\, k=1,2,3,4,$ are numbered in
consecutive order).

Estimate all integrals. Below we denote by
$c_i$ the positive constants, depending only on initial data.
Begin from $I_4.$ Before all we stress that according to
(2.9) $R'(\sigma)/R(\sigma)$ is bounded as
$\sigma\to\sigma_\infty,$ therefore $$|\frac{\partial\Psi^{-1}
(S,\sigma) }{\partial \sigma}|=|-\frac{\Psi^{-2}(S,\sigma)
\Psi(S,\sigma)
R'}{R}+ O(\Psi^{-1}(S,\sigma)
)|\le c_1\Psi^{-1}(S,\sigma).$$ Thus, $$I_4\le c_2(1+\|{\bf
W}{\bf \nabla}
_{\bf y}
s\|_\infty+ \|{\bf W}{\bf \nabla}_{\bf y}
\bar
S\|_\infty)E_p\le (c_2+c_3 R^{1/2} F_m+
c_4 R^{1/2} F_m^{1/2}
) E_p.$$

Further,
$$I_1=
-\int\limits_{R_n}D^\alpha P D^\alpha(({\bf W},{\nabla}_{\bf y})P+
\frac{\gamma-1}{2}(P+\bar P) {\rm div}_{\bf y}
 {\bf W}+
({\bf W},{\nabla}_{\bf y})\bar P+PQ_1(\sigma))d{\bf y}.$$
The integral from the member of higher order
$$D^\alpha P D^\alpha(({\bf W},{\nabla}_{\bf y})P+
\frac{\gamma-1}{2}P {\rm div}_{\bf y}
 {\bf W})$$
been integrated by parts can be reduced to
$$
-\int\limits_{R_n}\left\{\frac{1}{2}|D^\alpha P| {\rm div}_{\bf y}
 {\bf W}
+ D^\alpha P\left(({\bf W},{\nabla}_{\bf y}(D^\alpha P))-
D^\alpha({\bf W},{\bf \nabla}_{\bf y}P)\right)\right\} d{\bf y}-
$$
$$
-\frac{\gamma-1}{2}
\int\limits_{R_n}D^\alpha P\left\{\left ({\nabla}_{\bf y}P,
D^{\alpha}{\bf W})+(
P {\rm div}_{\bf y}
 D^{\alpha}{\bf W}-D^{\alpha}(P {\rm div}_{\bf y}
 {\bf W})\right )\right\}
d{\bf y}
:=I_{11}.
$$
Under sign of the integral there is the sum of members of the form
$\partial^\alpha U_1
\partial^l U_2 \partial^{\alpha+1-l} U_2,\, 0\le l \le \alpha $
(by
$U_j,\,j=1,2,$ we mean the corresponding components of
the solution.  At first let
$\alpha\ne 0,1.$ According to the Galiardo-Nirinberg inequality
we have
$$|\partial^i
U_j|_{p_i}\le C|DU_j|_\infty^{1-2/p_i}|D^\alpha U_j|_2^{2/p_i}$$
for $p_i=2\frac{\alpha-1}{i-1}.$
For $l\ne1,\alpha$ it is true that
$\frac{1}{p_l}+\frac{1}{p_{\alpha-l+1}}=\frac{1}{2}.$
From the H$\rm \ddot o$lder inequality
it follows that
$$\int\limits_{R_n} |\partial^\alpha U_1
\partial^l U_2  \partial^{\alpha+1-l} U_2| d{\bf y}\le $$ $$
\|D^\alpha U_1\|_2\|\partial^l U_2 \partial^{\alpha+1-l} U_2\|_2\le
$$
$$
\|D^\alpha U_1\|_2\|\partial^l U_2\|_{p_l}\| \partial^{\alpha+1-l}
U_2\|_{p_{\alpha-l+1}}\le
C \|D^\alpha U_1\|_2 \|D^\alpha U_2\|_2 \|D^\alpha U_2\|_\infty.
$$
For the other values of $\alpha$  and
$l$ the last inequality is evident.

Thus,
$$
I_{11}\le c\|D^p P\|_2(\|D^p P\|_2\|\nabla_ {\bf y}{\bf W}\|_\infty+
\|D^p {\bf W}\|_2\|\nabla_{\bf y} P\|_\infty)\le c
R^{1/2}F_m^{1/2}E_p.$$ (one can find the details of the inequalities
applying  in the book \cite{SerreBook}.)

Further, according to condition (2.7) of the Theorem
$$\left|\int\limits_{R_n}D^\alpha P D^\alpha ({\bf W}{\nabla}_{\bf y}
\bar P)d{\bf y}\right|\le c_5 R^{1/2}F_m,$$
$$\left|\int\limits_{R_n}D^\alpha P D^\alpha (PQ_1)d{\bf y}\right|\le
c_6 Q_1 R F_m\le c_7 F_m.$$

The integral $I_2$ can be estimated analogously:
$$|I_2|\le  c_8 R^{1/2}F_m^{1/2}E_q+ c_9 R^{1/2}F_m.$$

At last,
$$
I_3=
\int\limits_{R_n}
(D^\alpha)^*{\bf W}\Psi^{-1}
(S,\sigma)
\left \{({\bf W},{\nabla}_{\bf y}) D^\alpha{\bf W}
-D^\alpha(({\bf W}, {\nabla}_{\bf y}{\bf W})\right \}  d{\bf y}+
$$
$$
\frac{1}{2}\int\limits_{R_n}
\left\{(D^\alpha{\bf W})^*\Psi^{-1}
(S,\sigma)
D^\alpha{\bf W}) {\rm div}
_{\bf y}
{\bf W}
- \frac{{\bf W \nabla}_{\bf y}S}{\gamma}(D^\alpha{\bf W})^*
\Psi^{-1}(S,\sigma)
D^\alpha{\bf W}\right\} d{\bf y}+
$$
$$
\frac{\gamma-1}{2}
\int\limits_{R_n}
\left\{(D^\alpha{\bf W})^*\Psi^{-1}(S,\sigma)
D^\alpha\left (P\Psi(S,\sigma)
{\bf
\nabla}_{\bf y}P +\bar P \Psi (S,\sigma)
{\bf \nabla}_{\bf y} P- \right. \right.
$$
$$
\left.\left.
(P+\bar P) \Psi (S,\sigma){\bf
\nabla}_{\bf y} \bar P
+\bar P \Psi (\bar S,\sigma){\bf
\nabla}_{\bf y} \bar P
\right)
\right\}
d{\bf y}+ $$ $$
\int\limits_{R_n}
(D^\alpha{\bf W})^*\Psi^{-1}(S,\sigma)
D^\alpha\left ({\bf W}Q_2+G\right )
d{\bf y}.
$$

The first and second integrals can be estimated from above by
values

$$c_{10}R^{-1}\|D^p {\bf W}\|_2^2\|\nabla {\bf
W}\|_\infty+ c_{11}R^{-1}\|D^p {\bf W}\|_2^2\|{\bf W}\|_2) \le c_{12}
R^{1/2}F_m^{1/2}E_p,$$

and the third one by the value
$$
c_{13}\|D^p {\bf W}\|_2
(\|D^p P\|_2+\|D^p s\|_2+
\|D^p \bar P\|_2+\|D^{p-1}{\bf\nabla}_{\bf y}\bar S\|_2
)(1+
\|P\|_\infty+ \|{\bf\nabla}_{\bf y} P\|_\infty+$$$$
\|{\bf\nabla}_{\bf y} s+\|_\infty
\|\bar P\|_\infty+ \|{\bf\nabla}_{\bf y}\bar P\|_\infty+
\|{\bf\nabla}_{\bf y}\bar S\|_\infty
)^{p+1}
\le $$
$$c_{13} R^{1/2}(c_{14}+F_m^{(p+1)/2})E_p+
c_{15} R^{1/2}(c_{16}+F_m^{(p+1)/2})E_p^{1/2},$$
and, at last, the fourth one taking into account conditions (2.8)
by the value
$$c_{17} E_p.$$
In the last case note that
$$G=\lambda^{-1}(t)A(t){\bf f}_{\bf V} (t,{\bf x},\Pi,\bar {\bf
V}+\theta {\bf v}, S)A^{-1}(t),\,\theta\in(0,1),$$ $$
\int\limits_{R_n}
(D^\alpha{\bf W})^*\Psi^{-1}(S,\sigma)
D^\alpha(U_\phi(t){\bf W})
d{\bf y}=0,
$$
as $({\bf
w},U_\phi(t){\bf w})=0$ for any vector ${\bf w}\in
{\mathbb R}^n.$

Uniting all the estimates we get

$$E'_p\le c_{18}E_p+c_{19}R^{1/2}E_p + c_{20}R^{1/2}F_m^{1/2}E_p +
c_{21}R^{1/2}F_m E_p + $$
$$
c_{22}R^{1/2}F_m^{(p+1)/2}E_p+
c_{23}R^{1/2}F_m^{(p+1)/2}E_p^{1/2},
$$
whence it follows after summation over
$p$ that

$F_m'\le c_{24}F_m+c_{25}R^{1/2}(F_m^{1/2}+F_m+F_m^{3/2}+F_m^{(p+2)/2}+F_m^{(p+3)/2}).$

Set $\Lambda_m(\sigma)=\exp(-c_{24}\sigma)F_m.$
Then
$$\Lambda_m'\le c_{26}
(\Lambda_m^{1/2}+\Lambda_m+\Lambda_m^{3/2}+\Lambda_m^{(p+2)/2}
+\Lambda_m^{(p+3)/2})R^{1/2},\eqno(2.10)$$
where the constant $c_{26}$ depends on
$c_{24},\,c_{25},\,\sigma_\infty$.

We set
$$\Theta(\tilde g):=\int\limits_{\delta>0}^{\tilde g}
\frac{dg}{g^{1/2}+ g+g^{3/2}+g^{(p+2)/2}+g^{(p+3)/2}}.$$
The integral diverges in the zero, therefore
$\Theta(0)=-\infty.$

Integrating inequalities (2.10) over
$\sigma$ we obtain that $$ \Theta(\Lambda_m(\sigma))\le
\Theta(\Lambda_m(t_0))+ c_{26}\int\limits_{\sigma(t_0)}^\sigma
R^{1/2}(\tilde\sigma) d\tilde\sigma,$$
moreover, as follows from (2.8),
the integral in the right-hand side of the last inequality
converges as
$\sigma\to\sigma_\infty$ to the constant $C$, depending only on
the initial data.

Choosing $\Lambda_m(t_0),$ (and together with the value the
$H_m -$ norm of the initial data) sufficiently small,
one can obtain that $\Theta(\Lambda_m(t_0))+C$ is
later then $\Theta(+\infty),$ it signifies that $\Lambda_m(\sigma)$
and $F_m(\sigma)$ are bounded from above for all
$\sigma\in[0,\sigma_\infty) $ and $\sigma_*=\sigma_\infty$.  So,
Theorem 2.1 is proved.

\vskip1cm

{\sc Remark. } The solutions with the linear profile of
velocity, satisfying the theorem conditions, exist.
In the case described in \cite{Serre} $f=0,\,\gamma\le
1+\frac{2}{n},\,$ $A(t)=(E+A(0)t)^{-1}A(0),\, {\rm Spec} A(0)\notin
{\mathbb R_-},$ $\,{\rm tr} A(t)\sim \frac {n}{t},\,$$ {\rm det}
A(t)\sim \frac {1}{t^n},\,$ $Q_1=0,\,$$ Q_2\sim
$$2\left(\frac{1}{n}(({\rm tr} A^{-1})E-A^{-1}\right),\,$ as
$t\to\infty.$ Here
$\lambda=(1+t)^{-2},\,$$q=\frac{n(\gamma-1)}{4},$ $U_\phi(t)=0.$

\vskip1cm
                                            We shall get two corollaries from Theorem 2.1, basing on which
we can assert that the solutions which will be constructed in the next
section are interiors.

\begin{cor}
Let ${\bf f}={\bf 0}.$
If system
(1--3) has the globally smooth in time solution with the linear
profile of velocity, described in the statement of Theorem 2.1
and $A(t) \sim \frac{\delta}{t}E,\,t\to\infty,\,$
where $\delta  $ is a positive constant, then the solution is
interior.  \end{cor}

Proof.
Choose $\lambda(t)=t^{\delta+1},$
$q\ge\frac{3}{2}-\frac{n+1}{2(\delta+1)},$ $U_\phi(t)=0.$
In that way all conditions of Theorem 2.1 are satisfied.

\begin{cor}
Let ${\bf f}=L{\bf V},$
where $L $ is a matrix with the smooth coefficients such that
$A_0 L A_0^{-1}- U_1(t)=-\mu E,$
where  $A_0=\delta E+U_2$  is a matrix with
the positive determinant, the matrices
$U_1$ and $U_2$ are skew symmetric, $\delta $
is a constant, $\mu $ is a positive constant.
If system (1--3) has
the globally smooth in time solution with the linear profile of
velocity described in the statement of Theorem 2.1 and $A(t) \sim
\frac{1}{t}A_0,\,t\to\infty,\,$ then the solution is interior.
\end{cor}

Proof.
To verify condition (2.8) one has to show the boundedness of
the value
$$\lambda^{-1}(t)((\ln \lambda(t))'E+A(t)E-
A(t)LA^{-1}(t)-A^{-1}(t)A'(t)-U_\phi)\sim $$ $$\lambda^{-1}(t)((\ln
\lambda(t))'E+\frac{1}{t}((\delta+1) E+U_2)- U_1(t)+\mu
E-U_\phi),\,t\to\infty.$$
Choose $\lambda(t)=t^{\delta+1}\exp\{-\mu
t\}.$ Condition (2.7) is satisfied at
$q\ge\frac{3}{2},$ condition (2.9) can be verified
elementary. It needs to choose the matrix
$U_\phi(t)=\frac{1}{t}U_2-U_1(t).$

\section{Constructing of interior solutions}

Limit ourselves to the important case of
$n=2$ and set ${\bf f}= L{\bf v},$ where the matrix
$L=\left(\begin{array}{cc} -\mu&-l\\l&-\mu\end{array}\right),\, $ $
\mu$ is a nonnegative constant, $l$ is an arbitrary constant.
Thus, we in the simplest way  describe the Coriolis force and the
Rayleigh friction, for example, in the meteorological model under
neglecting the vertical processes and the Earth curvature.

It is easy to verify that on the smooth solutions to system
(1--3), so quickly vanishing at infinity to assure the convergency
of all integrals involved,
there are the following conservation quantities:
the mass $m=\int\limits_{{\mathbb R}^2} \rho\,d{\bf x}=const, $ for
$\mu=0$ the total energy $E=\int\limits_{{\mathbb R}^2}
\left( \frac{\rho |{\bf V}|^2)}{2}+\frac{p}{\gamma-1} \right) \,d{\bf
x}=E_k(t)+E_p(t)=const$ and the momentum of inertia
$\,J=\int\limits_{{\mathbb R}^2}  \rho(({\bf V}_\bot,{\bf r})+$
$\frac{l}{2}|{\bf x}|^2) \,d{\bf x}=const.  $ The conservation
laws are true also for the solutions to system (1--3) with the
density (and the pressure) quickly vanishing as $|{\bf x}| \to\infty$,
and the velocity components may even increase.

We consider another integral functionals,
characterizing the average properties of solutions,
namely $$G(t)=\frac{1}{2}\int
\limits_{{\mathbb R}^2}\rho|{\bf r}|^2\,d{\bf x},\quad
F_i(t)=\int\limits_{{\mathbb R}^2}({\bf V, X}_i)\rho\,d{\bf x},\quad
i=1,2,$$ where ${\bf X}_1={\bf r}=(x,y),\,{\bf X}_2={\bf
r}_\bot=(y,-x).$

   Note that $G(t)>0.$

On smooth solutions to system (1 -- 3)
the following relations take place
\cite{Roz98FAO}
$$ G'(t)=F_1(t), \eqno(3.1) $$ $$ F_2'(t)=l
F_1(t)-\mu F_2(t),\eqno(3.2) $$ $$ F_1'(t)=2(\gamma-1)E_p(t)+2 E_k(t)
- l F_2(t)-\mu F_1(t), \eqno(3.3) $$ $$E'(t)= -2\mu E_k(t). $$
Remark once more that at $\mu=0$ the quantity $E(t)$ of total
energy is constant.  The functions $E_p(t)$ and
$E_k(t)$ generally cannot be expressed through $G(t), F_1(t), F_2(t)$.

We try to choose the form of velocity field in such way
that from system (3.1
-- 3.3) one could find it directly.
Suppose, for example,

$$
{\bf V}= A(t){\bf r}, \eqno(3.4)
$$
where $A(t)$ is a (2 å 2) matrix , whose time-depending
coefficients we have to find.

1) We get the simplest result if we choose the matrix $A(t)$
such that

$$
{\bf V}= \alpha(t){\bf r}+\beta(t){\bf r_\bot},
$$
It is easy to see that in the case we have
$$
F_1(t)= 2\alpha(t) G(t), \quad  F_2(t)= 2\beta(t) G(t),
$$
$$
E_k(t)=(\alpha^2(t)+\beta^2(t)) G(t).
$$

It is not difficult to obtain
from the equation
$$ \partial_t p +({\bf
V},{\bf\nabla} p)+\gamma\,p\, div\,{\bf V}=0,  \eqno(3.5) $$
which is a corollary of system (1--3) and the state equation,
that
$$ E_p(t)= E_p(0)
G^{\gamma-1}(0) \frac{1}{G^{\gamma-1}(t)}.  $$

In that way, all functions involved in system (3.1
-- 3.3) are expressed through $G(t), \alpha(t),
\beta(t).$ For the sake of convenience we denote $G_1(t)=
1/G(t)$ and get the system of equations $$
G_1'(t)=-2\alpha(t) G_1(t), \eqno(3.5) $$ $$
\beta'(t)=\alpha(t)(l-2\beta (t)) - \mu \beta(t), \eqno(3.6) $$ $$
\alpha'(t)=-\alpha^2(t)+\beta^2(t)-l\beta(t)-\mu\alpha(t)+
(\gamma-1)E_p(0) G_1^{1-\gamma}(0) G_1^{\gamma}(t).  \eqno(3.7)
$$

The solutions of the system of ordinary differential equations
are smooth if they are can be prolonged up to all axis of time.
The non-prolongation is connected with the escape of the components
of solution to infinity during a finite time.

Note that for $\gamma > 1$ the quantities $\alpha(t)$ and
$\beta(t)$ are bounded, this follows immediately from the
expression for the total energy of the system:
$$
E(t)=(\alpha(t)^2+\beta(t)^2)G(t)+ E_p(0) G^{\gamma-1}(0)
\frac{1}{G^{\gamma-1}(t)}\le E(0), $$
therefore
$$
\alpha^2(t)+\beta^2(t)\le  {E(0)}{G_1(t)}- E_p(0) G_1^{1-\gamma}(0)
G_1^{\gamma}(t)<+\infty, $$
thus the quantities of density and pressure
remain bounded for all solution we want to construct.

At $\mu=0$ system (3.5--3.7) can be integrated explicitly:

$$
\alpha(G_1)=\pm\sqrt{KG_1^\gamma-C^2G_1^2+(E-lC)G_1-l^2/4},
$$

$$
\beta(t)=CG_1(t)+l/2, \eqno(3.8)
$$

$$ -\int_{G_1(0)}^{G_1(t)}
\frac{dG_1}{2G_1\alpha(G_1)} = t,
$$
with the constants
$C=\frac{2\beta(0)-l}{2G_1(0)},\quad
K=(\alpha^2(0)+C^2 G_1^2(0)-(E-lC)G_1(0)+l^2/4)/G_1^\gamma(0).$

Due to (3.8) one can reduce the number of equations in the system
and consider it on the phase plane
$(G_1, \alpha), G_1>0$.
There exists the unique singular point: at $l\ne0$ it is a center
situated on the axis $\alpha =0$, at $l=0$ it is a complex
equilibrium in the origin, the trajectories form the elliptic saddle
point.
In the last case the time of movement from any point of the phase
plane to the origin is infinite.
The equilibria, at $l\ne 0$ lying in the plane $G_1=0\,$, namely,
the points $\alpha=\beta=0 \,$ and $\alpha=0,\beta=l,$ are the centers
under consideration in the plane $(\alpha,\beta).$

At $\mu\ne 0$ one needs to integrate the system (3.5--3.7)
numerically.
However, the conclusion on the behaviour of the solution
components at infinity we can do analytically as well.
At $\mu>0$ the system has a unique stable equilibrium in the
origin (at $l=0$ it is a knot, at $l\ne 0$ it is a focus.)
Therefore we can find the solution as a formal asymptotic
series by the negative powers of $t$.
Limit ourselves to writing out the leading terms:

$$\alpha(t)\sim a_1 t^p,\,\beta(t)\sim b_1 t^r,\,G_1(t)\sim c_1 t^q,
\,a_1,b_1,c_1=const\ne 0.$$

From (3.5) we have immediately that
$q c_1 t^{q-1}=-2a_1 c_1 t^{p+q},$
therefore $p=-1,\, q=-2a_1.$

From (3.6) we obtain that
$r b_1 t^{r-1}\sim l a_1 t^p - \mu b_1 t^r -2 a_1 b_1 t^{p+r}.$
If $l\ne 0,\,\mu\ne 0,$ then $r=p=-1,$ moreover, $la_1=\mu b_1.$

From (3.7) we get analogously that
$pa_1 t^{p-1} \sim -a_1^2 t^{2p}+ b_1^2 t^{2r} - lb_1t^r-\mu a_1 t^p
+Kc_1^\gamma t^{\gamma q},$ where it is denoted for short
$K=(\gamma-1)E_p(0) G_1^{1-\gamma}(0),$ that is
$-a_1 t^{-2} \sim -a_1^2 t^{-2}+ b_1^2 t^{-2} - lb_1t^{-1}-\mu
a_1 t^{-1} +Kc_1^\gamma t^{\gamma q}.$

Let $\gamma q<-1,$ that is the last term does not contain the
senior degree.
Then the rest of members of senior order must be eliminated, and
$-lb_1=\mu a_1=\mu^2 b_1/l,$ therefore $l=\mu=0$ in spite of
the supposition.
Consequently,
$q=-1/\gamma.$

So, at $\mu\ne 0, \,l\ne 0$
$$\alpha(t)\sim \frac{1}{2\gamma}t^{-1},\,
\beta(t)\sim \frac{\mu}{2l\gamma}t^{-1},\,
G_1(t)\sim \left( \frac{\mu^2+l^2}{2K\gamma\mu}\right)^{1/\gamma}
t^{-1/\gamma},\,t\to\infty.$$

At $l=0,$ due to the existence of the integral
$\beta(t)=CG_1(t) \exp\{-\mu t\},$
system (3.5 -- 3.7) can be reduced to the nonautonomous
system of two equations:
$$ G_1'(t)=-2\alpha(t) G_1(t), $$ $$
\alpha'(t)=-\alpha^2(t)-\mu\alpha(t)+ C^2 \exp \{-2\mu t\}G_1^2(t)+ K
G_1^{\gamma}(t).  \eqno(3.9) $$
From (3.9) we have that $pa_1 t^{p-1}
\sim -a_1^2 t^{2p}-\mu a_1 t^p+C^2c_1^2 t^{2q}\exp\{-2\mu t\}
+Kc_1^\gamma t^{\gamma q},$ where, remember, $p=-1,\,
q=-2a_1.$ To compensate the senior term $-\mu a_1
t^{-1}$ it is necessary that the equality
$ q=-1/\gamma$ holds. Hence, we find the coefficient
$c_1.$

Thus, at $\mu\ne 0, \,l=0$
$$\alpha(t)\sim \frac{1}{2\gamma}t^{-1},\,
\beta(t)\sim C\left(\frac{\mu}{2K\gamma}\right)^{1/\gamma}t^{-1/\gamma}\exp\{-\mu
t\},\, G_1\sim \left(\frac{\mu}{2K\gamma}\right)^{1/\gamma}
t^{-1/\gamma},\,t\to\infty.$$

At $\mu=l=0,$ as follows from the explicit formulas,
$$
\alpha(t)\sim\left(\sqrt{G_1(0)/E}+t\right)^{-1},\, G_1(t)\sim
\alpha^2(t),\, \beta(t)= C G_1(t),\,t\to\infty.$$

Thus, as follows from Corollaries 2.1 and
2.2, in the cases $\mu=l=0$ and $\mu>0 $ we have constructed
the velocity field for the interior solution.
More precisely, at $\mu=l=0$ $A(t)\sim \frac{1}{t}E,\,
t\to\infty,$ Corollary 2.1 can be applied with $\delta=1.$
At $\mu>0 $ one can apply Corollary 2.2. In the case
$\mu>0,\,l=0$ we have $U_1=U_2=0,$ $\delta=\frac{1}{2\gamma}$.
If $\mu>0,\,l\ne 0,$ then $\delta=\frac{1}{2\gamma},$
$U_1=\left( \begin{array}{cc} 0 & l\\ -l & 0\end{array}\right),\,$ $
U_2=\frac{\mu}{2\gamma l}\left( \begin{array}{cc} 0 & 1\\ -1 &
0\end{array}\right).  $

As soon as $\alpha(t)$ and $\beta(t)$ are found, the rest of the
components can be found elementary
$$ \rho(t,|r|,
\phi)=\exp(-2\int\limits_0^t \alpha(\tau) d\tau)
\rho_0(|r|\exp(-\int\limits_0^t \alpha(\tau) d\tau),
\phi+\int\limits_0^t \beta(\tau) d\tau), $$ $$ S(t,|r|, \phi)=
S_0(|r|\exp(-\int\limits_0^t \alpha(\tau) d\tau),
\phi+\int\limits_0^t \beta(\tau) d\tau).  $$
$$
p(t,|r|, \phi)=\exp(-2\gamma\int\limits_0^t \alpha(\tau) d\tau)
p_0(|r|\exp(-\int\limits_0^t \alpha(\tau) d\tau),
\phi+\int\limits_0^t \beta(\tau) d\tau). $$
From (2) and (3.7) we get that on the classical solution of the
initial system
must be satisfied the relation
$${\bf \nabla} p= -
(\gamma-1)G_1^{1-\gamma}(0)E_p(0)G_1^\gamma(t)\rho {\bf r}.$$
Hence it follows
that the components of the initial data
$\rho_0$ ¨ $p_0$ must be axisymmetric and compatible,
i.e. connected as follows:
$${\bf
\nabla} p_0= - (\gamma-1)G_1(0)E_p(0)\rho_0 {\bf r}.$$

However, in spite of the components of the pressure and the
density must vanish as
$|{\bf x}|\to\infty,$ the entropy must even increase
(remember that the conditions to Theorem 2.1
require  only the boundedness of the entropy gradient).
However, one can choose $$p_0=\frac{1}{(1+|{\bf
r}|^2)^a},\,a=const>3,\quad \rho_0=
\frac{2a}{(\gamma-1)G_1(0)E_p(0)}\frac{1}{(1+|{\bf r}|^2)^{a+1}},$$
then $S_0=const+(a(\gamma-1)+\gamma)\ln(1+|{\bf r}|^2).$

{\sc Remark.} In the case $G(0)\ne 0$ the
density cannot be compactly supported
(in contrast $G(0)=0$).  Really, since
$p=\pi^\frac{2\gamma}{\gamma-1},\, \rho=
\pi^\frac{2}{\gamma-1}e^{-\frac{S}{\gamma}},$ then due to the
compatibility conditions
$\pi\nabla\pi\sim const\cdot e^{-\frac{S}{\gamma}}, \,|{\bf x}|\to
c-0,$ where $c $ is a point of the support of $\pi.$
Therefore, for $C^1 -$ smooth $\pi$ it occurs that $S\to
+\infty, \,|{\bf x}|\to c-0,$ and we cannot choose
any smooth initial data.  It is interesting that
if one requires only $C^0 -$ smoothness of
$\pi\,(\pi\sim const\cdot (c-|{\bf x}|)^{1/2},\,|{\bf x}|\to c-0, \,
\pi=0,\,|{\bf x}|\ge c>-0),$ the condition may be fulfilled.
Moreover, $\rho$ and $p$ will be of the $C^1$ class of smoothness,
however, neither the theorem on the local in time existence of the
smooth solution, no Theorem 2.1 can be applied.

\vskip1cm

2) To consider the velocity field
(3.4) with the matrix

$$A(t)=\left(
\begin{array}{cc} a(t) & b(t)\\ c(t) &
d(t)\end{array}\right),
$$
involve the functions
$$ G_x(t)=\frac{1}{2}\int\limits_{{\mathbb R}^2} \rho x^2
d{\bf x},\quad
G_y(t)=\frac{1}{2}\int\limits_{{\mathbb R}^2} \rho y^2 d{\bf x},
\quad G_{xy}(t)=\frac{1}{2}\int\limits_{{\mathbb R}^2} \rho xy
d{\bf x}.  $$
Then involve the auxiliary variables
$$G_1(t)={G_x(t)}{\Delta^{-(\gamma+1)/2}(t)},\,
G_2(t)={G_y(t)}{\Delta^{-(\gamma+1)/2}(t)},\,
G_3(t)={G_{xy}}{\Delta^{-(\gamma+1)/2}(t)},\eqno(3.10)$$
where $\Delta(t)=G_x(t)G_y(t)-G_{xy}^2(t)$ is a positive
function on solutions to system (1--3).  Note that the behaviour of
$\Delta(t)$ is governed by the equation
$$
\Delta'(t)=2(a(t)+d(t))\Delta(t),\eqno(3.11)$$
and a potential energy
$E_p(t)$ is connected with $\Delta(t)$ as follows:
$$E_p(t)=E_p(0)\Delta^{(\gamma-1)/2}(0)\Delta^{(-\gamma+1)/2}(t).$$

To find $G_1,\,G_2,\,G_3$, and also the elements of the matrix
$A(t)$ we get the system of equations
$$
G_1'(t)=((1-\gamma)a(t)-(1+\gamma)d(t))G_1(t)+2b(t)G_3(t), $$ $$
G_2'(t)=((1-\gamma)d(t)-(1+\gamma)a(t))G_2(t)+2c(t)G_3(t),
$$
$$
G_3'(t)=c(t)G_1(t)+b(t)G_2(t)-\gamma(a(t)+d(t))G_3(t),
$$
$$
a'(t)=-a^2(t)-b(t)c(t)+lc(t)-\mu a(t)+K_1G_2(t),\eqno(3.12)
$$
$$
b'(t)=-b(t)(a(t)+d(t))+ld(t)-\mu b(t)-K_1 G_3(t), $$
$$
c'(t)=-c(t)(a(t)+d(t))-la(t)-\mu c(t)-K_1G_3(t), $$ $$
d'(t)=-d^2(t)-b(t)c(t)+lb(t)-\mu d(t)+K_1G_1(t), $$
with
$K_1=\frac{\gamma-1}{2}
(E_p(0)\Delta^{(\gamma-1)/2}(0))$.

\newtheorem{utv}{Proposition 3.}
\begin{utv}
In the case $\mu=l=0$ any solution to system
(1--3) with a linear profile of velocity,
satisfying the condition  )
of Theorem 2.1,
is interior.  \end{utv}

Proof.
Show that at $\mu=l=0$ there is the asymptotics
$A(t)\sim \frac{\delta}{t}E,\,t\to\infty,$ that is
according to Corollary 2.1 the corresponding solution
with a linear profile of velocity
is interior.

Go to the new variables
$a_1(t)=a(t)-d(t), b_1(t)= b(t)+c(t), c_1(t)=b(t)-c(t),
d_1(t)=a(t)+d(t), G_4(t)=G_1(t)+G_2(t), G_5(t)=G_1(t)-G_2(t).$
In  variables
$ a_1,b_1,c_1,d_1,G_3,G_4,G_5 $ system (3.12) has a form:
$$ a_1'(t)=-a_1(t)d_1(t)-K_1G_5(t),\eqno(3.13) $$ $$
b_1'(t)=-b_1(t)d_1(t)-2K_1G_3(t),\eqno(3.14)
$$
$$
c_1'(t)=-c_1(t)d_1(t),\eqno(3.15)
$$
$$
d_1'(t)=-\frac{1}{2}(a_1^2(t)+d_1^2(t))-\frac{1}{2}(b_1^2(t)-c_1^2(t))+
K_1G_4(t),\eqno(3.16)
$$
$$
G_3'(t)=-\gamma d_1(t)G_3(t)+\frac{1}{2}
b_1(t)G_4(t)-\frac{1}{2}c_1(t)G_5(t),\eqno(3.17) $$
$$
G_4'(t)=-\gamma
d_1(t)G_4(t)
+a_1(t)G_5(t)+2b_1(t)G_3(t),\eqno(3.18) $$
$$
G_5'(t)=-\gamma
d_1(t)G_5(t)+a_1(t))G_4(t)-2c_1(t)G_3(t).\eqno(3.19) $$

Let as $ t\to\infty $ the asymptotics of functions involved in the
system be the following:  $
a_1(t)\sim L_1 t^{l_1}, $ $b_1(t)\sim L_2 t^{l_2}, c_1(t)\sim L_3
t^{l_3}, $$\, d_1(t)\sim L_4 t^{l_4}, $$\,G_3(t)\sim N t^{q},\,
G_4(t)\sim M_1 t^{p_1}, $ $\,G_5(t)\sim M_2 t^{p_2},  $ where $L_i,
\,i=1,2,3,\, M_j,\, j=1,2, \, N $ are some constants not
equal to zero.  Note that
$p_2\le p_1$ and in virtue of $\Delta>0$
the estimate $q\le p_1$ holds.

From (3.15) we get immediately that
$$l_3 L_3 t^{l_3-1}=-L_3 L_4 t^{l_3+l_4},$$
hence
$l_4=-1,\, L_4=-l_3.$

From (3.13) we have taking it into account
$$l_1 L_1 t^{l_1-1}=-L_1 L_4 t^{l_1-1}-2K_1 M_2 t^{p_2}.$$
The following variants are possible:

$p_2<l_1-1, $ hence $l_1=-L_4,$ i.e. $l_1=l_3,$

$p_2=l_1-1, $ hence
$$l_1=-L_4-\frac{K_1M_2}{L_1}=l_3-\frac{K_1M_2}{L_1}.\eqno(3.20)$$

From (3.14) we have analogously
$$l_2 L_2 t^{l_2-1}=-L_2 L_4 t^{l_2-1}-2K_1 N t^q.$$
There are the variants:

$q<l_2-1, $ hence $l_2=-L_4,$ i.e. $l_2=l_3,$

$q=l_2-1, $ hence
$$l_2=-L_4-\frac{2K_1N}{L_2}=l_3-\frac{2K_1N}{L_2}.\eqno(3.21)$$

From (3.16) we get
$$-L_4 t^{-2}=-\frac{1}{2}(L_1^2 t^{2l_1}+L_4^2 t^{-2})-
\frac{1}{2}(L_2^2 t^{2l_2}-L_3^2 t^{2l_3})+
K_1M_1t^{p_1}.\eqno(3.22)$$

If $l_i<-1,\,i=1,2,3,$ then at $t\to\infty$
$$a(t)\sim\frac{1}{2}(L_1t^{l_1}+L_4 t^{-1}),\quad
b(t)\sim\frac{1}{2}(L_2t^{l_2}+L_3 t^{l_3}),$$
$$c(t)\sim\frac{1}{2}(L_2t^{l_2}-L_3 t^{l_3}),\quad
d(t)\sim\frac{1}{2}(L_4t^{-1}-L_1 t^{l_1}),$$
that is $A(t)\sim \frac{L_4}{2t}E,\, L_4=const>0.$

We shall show below that there are not others situations.

From (3.22) it follows the impossibility of the situation
$l_3=-1, l_1<-1,l_2<-1,p_1\le -2.$ Really, in the case
$l_3=-L_4=-1-\sqrt{1+L_3^2}<-1$ or
$l_3=-L_4=-1-\sqrt{1+L_3^2+2K_1M_1}<-1,$
it contradicts to the supposition.

Now let $l_3>-1.$ Then from (3.22)
$l_1=l_3$ or (and) $l_2=l_3.$
Suppose, for example, that $l_1=l_3.$
Then $p_2<l_1-1,$ from (3.19) it follows that
$$p_2M_2 t^{p_2-1}=-\gamma L_4 M_2 t^{p_2-1} + L_1M_1t^{p_1+l_1}
+2L_2Nt^{q+l_2},
\eqno(3.23)
$$
and therefore
$p_1+l_1=q+l_2,\,l_1\le\l_2,$
i.e. $l_2>-1.$
Besides,
$$L_1M_1=-2L_2N.\eqno(3.24)$$

From (3.17) we have
$$qN t^{q-1}=-\gamma L_4 N t^{q-1} + \frac{1}{2}L_2M_1t^{p_1+l_2}
-\frac{1}{2}L_3M_2t^{p_2+l_3}.
\eqno(3.25)
$$
As soon as $ q-1\le p_1-1, $ and  $ p_1-1 < p_1+l_2,$
then
$p_1+l_2=p_2+l_3\le p_1+l_3, $
and therefore
$l_2\le l_3,$
and as
$l_1=l_3, $ then $l_1=l_2=l_3$ and $p_1=q.$
Besides,
$$L_2 M_1=L_3M_2.\eqno(3.26)$$

Further, from (3.18) we obtain
$$p_1M_1 t^{p_1-1}=-\gamma L_4 M_1 t^{p_1-1} + L_1M_2t^{p_2+l_1}
-2L_3Nt^{q+l_3}.
\eqno(3.27)
$$
As $ q+l_3> p_1-1, $ then
$p_2+l_1=q+l_3,\, p_2=q, $
$$L_2 M_1=2L_3N.\eqno(3.28)$$

From (3.24) and (3.28) we have
$\frac{M_1}{M_2}=-\frac{L_2}{L_3},$ and from (3.26) we
obtain  $\frac{M_1}{M_2}=\frac{L_2}{L_3},$ i.e.
$L_2^2=-L_3^2,\, L_2=L_3=0$ in spite of the supposition.

Suppose now that $l_2=l_3.$ Then
$q<l_2-1,\,p_2\le l_1-1.$ From (3.17) we have that $q-1<l_2+p_1,
$ since $l_2>-1,$ and therefore $l_2+p_1=l_3+p_2,\,
p_1=p_2,\,L_2M_1=L_3M_2.$

Taking this into account from (3.23)  we have that
$p_1+l_1=q+l_2\le p_1+l_2,\,l_1\le l_2,\, L_2M_1=-2L_2N.$

From (3.27) we get
$p_2+l_1=q+l_3\le p_1+l_2,\,p_2=q,\,l_1=l_2,\, L_1M_2=-2L_3N.$
In that way,  we obtain the contradiction
analogous to the previous one.

Now let $l_3\le-1,$   $l_1>-1 $ and (or)   $l_2>-1.$
Then from (3.22) it follows that
$p_1=2l_1 $ and (or)   $p_1=2l_2\, (p_1>-2). $
For example, if
$p_1=2l_1, $ then from (3.24) we get
$p_2-1<p_1+l_1,\, p_1+l_1=q+l_2\le p_1+l_2,\,l_1\le l_2.$
From (3.25)
$p_2+l_3\le p_1-1,\,l_2+p_1\ge l_1+p_1>p_1-1,$
therefore
$q-1=l_2+p_1.$
But
$q-1\le p_1-1,$
therefore $l_2\le-1,$ and $l_2\le-1$ in spite of the
supposition.

If $p_1=2l_2, $ then from (3.25) we get
$l_3+p_2\le p_1-1,$ therefore $p_1+l_2>=p_2+l_3,\,
q-1=p_1+l_2,\,l_1\le l_2.$ But $q-1\le p_1-1, $
therefore $l_2\le-1$ in spite of the supposition.

So, it remains the unique possibility:
$l_3\le -1,\,l_1\le -1 $ and (or) $l_2\le -1.$
In the case, as follows from (3.22)
$$l_3=-1\pm\sqrt{1-(\delta_1 L_1^2+\delta_2 L_2^2 -\delta_3 L_3^2
- 2\delta_4 K_1 M_1)},\eqno(3.29)$$
where $  \delta_i=1,$ if $
l_i=-1, $ and $ \delta_i=0 $ otherwise, $
i=1,2,3, \, $  $ \delta_4=1, $ if $p_2=-2 $  and $ \delta_4=0$
otherwise. Hence
$$\delta_1 L_1^2+\delta_2 L_2^2 -\delta_3 L_3^2
- 2\delta_4 K_1 M_1\le 1, \eqno(3.30)$$
if inequality (3.30) is strict, then $l_3<-1.$

We consider this case, that is
$l_3<-1,\,l_1=-1$ and (or)  $l_2=-1.$

If $l_1=-1\ne\l_3,$ then $p_2=l_1-1=-2.$
If $l_2=-1\ne\l_3,$ then
$q=l_2-1=-2.$  In that way, in any of these cases we have
$ p_1=-2.$

At $l_1=-1,\,l_2<-1,\,l_3<-1,$
from (3.23)
we obtain
$p_2=-\gamma L_4+\frac{L_1M_1}{M_2},$
$$\frac{L_1M_1}{M_2}=-(2+\gamma l_3),\eqno(3.31)$$
from (3.27)  we get
$p_1=-\gamma L_4+\frac{L_1M_2}{M_1},$
$$\frac{L_1M_2}{M_1}=-(2+\gamma l_3),\eqno(3.32),$$
and from (3.21)
$$
\frac{K_1M_2}{L_1}=l_3+1. \eqno(3.33)$$
From (3.31), (3.33) after excluding  $L_1$ and $M_2$
we obtain
$$K_1M_1=-(2+\gamma l_3)(l_3+1).\eqno(3.34)$$

As $K_1$ and $M_1$ are positive, and $l_3\le -1,$
then
$$l_3>-\frac{2}{\gamma}.\eqno(3.35)$$

At $\gamma>2$ it goes already to the contradiction.

Further, from (3.32) and (3.33) we have
$\frac{L_1^2}{K_1M_1}=-\frac{2+\gamma l_3}{l_3+1},$
this fact together with (3.34) give
$$L_1^2= (2+\gamma l_3)^2. \eqno(3.36)$$

From (3.31), (3.34), (3.35) and (3.36) we obtain the
inequality $$-\frac{2}{\gamma}<-1-\sqrt{1-
(2+\gamma l_3)^2 -2 (2+\gamma l_3)(l_3+1)},\eqno(3.37)$$
which cannot be true at $\gamma>1.$

       If $l_1<-1,\,l_2=-1,\,l_3<-1,$ then from (3.23)
we get
$p_2=-\gamma L_4+\frac{2L_2N}{M_2},$
$$
\frac{2L_2N}{M_2}=-2-\gamma l_3, \eqno(3.38)
$$
from (3.25)
$q=-\gamma L_4+\frac{L_2M_1}{2N},$
$$
\frac{L_2M_1}{2N}=-2-\gamma l_3, \eqno(3.39)
$$
from (3.27)
$$p_2=-\gamma L_4=\gamma l_3.$$
In that way, $l_3=-\frac{2}{\gamma},$
this fact together with (3.38) or (3.39) contradicts to
the fact that $L_2, \,N$ and $M_1$ are not equal to zero.

    If $l_1<-1,\,l_2=-1,\,l_3=-1,$ then from (3.23),(3.25) and
(3.27) we get correspondingly
$$\frac{2L_2N}{M_2}+\frac{L_1M_1}{M_2}=-2-\gamma l_3,\eqno(3.40)$$
$$
\frac{L_2M_1}{2N}=-2-\gamma l_3, \eqno(3.41)
$$
$$
\frac{L_1M_2}{M_1}=-2-\gamma l_3, \eqno(3.42)
$$
from (3.20), (3.21)
$$
\frac{K_1M_2}{L_1}=\frac{2K_1N}{L_2}=l_3+1. \eqno(3.43)$$
In that way, multiplying (3.40) by (3.42), taking into
account (3.43) we get
$$ (2+\gamma
l_3)^2=L_1^2+\frac{M_2L_2^2}{M_1},$$ from (3.41), (3.42), (3.43)
$$M_1^2=M_2^2,\eqno(3.42)$$
that is $$ (2+\gamma l_3)^2=L_1^2\pm
L_2^2.\eqno(3.45) $$
As above, from (3.31), (3.34), (3.35) and (3.45)
we get the inequality (3.37) which cannot hold at
$\gamma>1.$

It remains the unique possibility
$l_1=l_2=l_3=-1.$
In the case $p_2<-2,\,q<-2$ and, as follows from (3.23),
$p_1=p_2$ or $p_1=q,$ that is $p_1<-2.$ Therefore, as follows
from (3.29)
$$ L_1^2+L_2^2-L_3^2=1.\eqno(3.46)$$

If $p_1=p_2$, then $|M_2|<M_1.$
If  $p_1=p_2,\,q<p_1,$ then from (3.23), (3.27)
we obtain that
$p_2+\gamma=\frac{L_1M_1}{M_2},\,
p_1+\gamma=\frac{L_1M_2}{M_1},$
hence $M_1^2=M_2^2,$ it goes to the contradiction.

If $p_1=q,\,p_2<p_1,$
then
from (3.23), (3.25), (3.27)
we get
       $$L_1M_1=-2L_2N,\,
\lambda=
\frac{L_2M_1}{2N}=-\frac{2L_3N}{M_1},\eqno(3.47)$$
where
$\lambda= q+\gamma
=p_1+\gamma.$
It follows, in particular, that
$L_2^2=-\lambda L_1,\,L_1L_3=-\lambda L_2,
\,\lambda^2=-L_2L_3.$
Besides, if $G_1(t)\sim N_1
t^{q_1},\,G_2(t)\sim N_2 t^{q_2},\,t\to\infty, $ where $N_1, N_2
$ are some positive constants,
$q_1=q_2=p_1=q,\, N_1=N_2,\, M_1=2N_1, $ and $N_1^2\ge N^2 $ in
virtue of $G_1(t)G_2(t)-G_3^2(t)>0.$ In such way, from
(3.47) we get that $L_2^2\ge \lambda^2,\, \lambda^2 \ge
L_3^2,\,L_1^2\ge L_2^2,$ and taking into account (3.46),
$L_1^2\le 1,\, L_2^2\le 1,\,\lambda^2\le 1.$ That is if
$\lambda>1\,(\gamma>3),$ the conditions mentioned in  the
paragraph cannot hold together.

At last, consider the case $p_1=p_2=q.$ Then from
(3.23), (3.25), (3.27) we have
$\lambda=\frac{L_1M_1}{M_2}+\frac{2L_2N}{M_2}=\frac{L_2M_1}{2N}-
\frac{L_3M_2}{2N}=\frac{L_1M_2}{M_1}-\frac{2L_3N}{M_1}, $
that is the system of linear homogeneous equations with respect to
the variables $M_1,\,M_2,\,N$ $$ L_1M_1-\lambda M_2+2L_2N=0,$$ $$
L_2M_1-L_3M_2-2\lambda N=0,$$ $$ \lambda M_1-L_1M_2+2L_3N=0.$$
For the existence of its nontrivial solution the determinant of
the system must be equal to zero, i.e.
$$
\lambda L_1^2-2\lambda L_2L_3+L_1L_2^2+L_1L_3^2-\lambda^3=0.$$
Taking into account (3.46), involve the function with respect
to the variables $L_2$ and
$L_3$, where $\lambda$ plays the role of parameter,
namely
$$\Psi_{\lambda}(L_2,L_3)=\lambda(1-L_2^2+L_3^2)-2\lambda L_2 L_3 +
(L_2^2+L_3^2)\sqrt{1-L_2^2+L_3^2}-\lambda^3.$$
By the standard methods one can show
that the function is not equal to zero at
$\lambda>1,$ i.e. at $\gamma\ge 3$
(remember, that $p_1<-2$).

So, it remains to investigate the case
$1<\gamma\le 3.$

Consider equation (3.11), which can be written as
$$ \Delta'(t)=2d_1(t)\Delta(t).$$
If we suppose that
$ \Delta(t)\sim const\cdot t^m,\,t\to\infty,\,m $ is a
constant, then $ m=2L_4=2.$ From (3.11) we can get
$E_p(t)\sim const\cdot t^{1-\gamma}(\to 0), \,t\to\infty.$
Therefore,
if we denote $E$ the quantity of the total energy of the
system, then $E_k(t)=E-E_p(t)\sim
E(1-C_1 t^{1-\gamma}),$ here and further $C_i $
are
some positive constants.  But
$E_k(t)\sim C_2 G_i(t) t^{-2}\Delta^{\frac{\gamma+1}{2}}\sim C_3 G_k
t^{\gamma-1},$ where $G_k  $ is at least one of functions
$G_1,\,G_2$ or $G_3.$ That is $G_k\sim
C_4(1-C_1 t^{1-\gamma})t^{1-\gamma}\sim C_5 t^{1-\gamma}.$ But then
$1-\gamma\le p_1<-2,\, \gamma>3, $ and $\gamma$
do not belong to the interval under consideration.

So the proof of the proposition is over.

\medskip

{\sc Remark 3.1.} One can  show more shortly, that in the physical
case $1<\gamma \le 2$ the situation $l_1=l_2=l_3=-1$ if impossible.
Taking into account (3.10) we have $G_1G_2-G_3^2=\Delta^{-\gamma}\sim
const\cdot t^{-2\gamma},\, t\to\infty.$ But the degree of the leading
term of the expression $G_1G_2-G_3^2$ is not greater then
$2p_1,$ therefore $p_1\ge-\gamma,$ as $p_1<-2,$ then $\gamma>2.$
\medskip

{\sc Remark 3.2.} Actually $l_3=-L_4=-2,\,
G_1G_2-G_3^2\sim const\cdot t^{-4\gamma},\,t\to\infty.$

{\sc Remark 3.3.} In the case $\mu>0,\,l=0$ the asymptotic
$A(t)\sim \frac{\delta}{t}E,\,t\to\infty$ holds as well, that is
according to Corollary 2.1 the corresponding solution
with a linear profile of velocity is interior.  If $\mu>0,\,l\ne
0$ it occurs that $A(t)\sim A_0 \frac{1}{t},$ but the matrix
$A_0$ under arbitrary initial conditions $A(0)$ has not to be of
form $\delta E+ U_2$ (see the denotations in the statement
of Corollary 2.2), i.e. the Corollary cannot be applied in the last
case.

\medskip

One can write out the conditions of compatibility of initial
data of the density and the pressure for the case of solution
with arbitrary linear profile of velocity.

\end{document}